\documentclass{amsart}
\usepackage{enumerate}
\input{diagrams.tex}
\diagramstyle[scriptlabels,small]

\def \C{{\mathbb C}}

\def \Q{{\mathbb Q}}
\def \S{{\bf S}}
\def \V{{\mathbb V}}
\def \Z{{\mathbb Z}}

\def\Proj{\mathop{\rm Proj}\nolimits}

\def\Spec{\mathop{\rm Spec}\nolimits}

\def\rk{\mathop{\rm rk}\nolimits}
\def\rank{\mathop{\rm rank}\nolimits}
\def\depth{\mathop{\rm depth}\nolimits}

\def\Hom{\mathop{\rm Hom}\nolimits}

\def\trdeg{\mathop{\rm trdeg}\nolimits}
\def\cHom{\mathop{\mathcal Hom}\nolimits}

\def\Pf{\mathop{\rm Pf}\nolimits}

\def\lto{\longrightarrow}

\newcommand{\ul}{\underline}
\newcommand{\ol}{\overline}

\newcommand{\cH}{{\mathcal H}}
\newcommand{\cI}{{\mathcal I}}

\newcommand{\cL}{{\mathcal L}}
\newcommand{\cM}{{\mathcal M}}
\newcommand{\cN}{{\mathcal N}}
\newcommand{\cO}{{\mathcal O}}

\newcommand{\cR}{{\mathcal R}}

\newcommand{\fa}{{\mathfrak a}}

\newcommand{\fm}{{\mathfrak m}}

\newcommand{\fq}{{\mathfrak q}}
\newcommand{\fp}{{\mathfrak p}}
\newcommand{\fM}{{\mathfrak M}}

\newcommand{\fP}{{\mathfrak P}}

\theoremstyle{definition}
\newtheorem{dfn}{Definition}[section]

\theoremstyle{plain}
\newtheorem{pro}[dfn]{Proposition}
\newtheorem{thm}[dfn]{Theorem}
\newtheorem{lem}[dfn]{Lemma}
\newtheorem{cor}[dfn]{Corollary}

\theoremstyle{remark}
\newtheorem{rem}[dfn]{Remark}
\newtheorem{sit}[dfn]{}
\newtheorem{exa}[dfn]{Example}
\newtheorem{que}[dfn]{Question}

\def\bque{\begin{que}}
\def\eque{\end{que}}
\def\bpro{\begin{pro}}
\def\epro{\end{pro}}
\def\blem{\begin{lem}}
\def\elem{\end{lem}}
\def\brem{\begin{rem}}
\def\erem{\end{rem}}
\def\bsit{\begin{sit}}
\def\esit{\end{sit}}
\def\bexa{\begin{exa}}
\def\eexa{\end{exa}}
\def\bcor{\begin{cor}}
\def\ecor{\end{cor}}
\def\bthm{\begin{thm}}
\def\ethm{\end{thm}}
\def\bdfn{\begin{dfn}}
\def\edfn{\end{dfn}}
\def\bproof{\begin{proof}}
\def\eproof{\end{proof}}
\def\bnum{\begin{enumerate}}
\def\enum{\end{enumerate}}

\def\bdi{\begin{diagram}}
\def\edi{\end{diagram}}

\def\ba{\begin{array}}
\def\ea{\end{array}}

\def\be{\begin{equation}}
\def\ee{\end{equation}}

\def\coker{\mathop{\rm coker}\nolimits}
\def\Hom{\mathop{\rm Hom}\nolimits}

\begin{document}

\title[Codimension and connectedness]{Codimension and
connectedness of degeneracy loci over local rings}
\author{Hubert Flenner}
\thanks{This paper originated during visits of the authors at
their home universities. They are grateful for the hospitality and
the financial support by the DFG Schwerpunkt ``Globale Methoden
der komplexen Analysis". The second author was supported in part
by the NSF}
\address{Fakult\"at f\"ur Mathematik der Ruhr-Universit\"at,
Universit\"atsstr.\ 150, Geb.\ NA 2/72, 44780 Bochum, Germany}
\email{hubert.flenner@rub.de}
\author{Bernd Ulrich}
\address{Department of Mathematics, Purdue University,
West Lafayette, IN 47907-1395, USA}
\email{ulrich@math.purdue.edu}


\begin{abstract}
We deduce results on the dimension and connectedness of degeneracy
loci of maps of finite modules $f:M\to N$ over a local noetherian
ring $(A,\fm)$. We show for instance that the expected
determinantal bounds on the dimension of the t-$th$ degeneracy
locus of $f$ hold if $f\in \fm\Hom (M,N)$, and that this
degeneracy locus is connected in the expected dimension provided
$\hat A$ is a domain.
\end{abstract}

\maketitle

\section{Introduction}

In this paper we will investigate the dimension and connectedness
properties of degeneracy loci of homomorphisms $f:M\to N$ of
finite modules over a noetherian local ring $(A,\fm)$. For $t\ge
0$ we consider the $(t+1)$-$th$ determinantal ideal $I_{t+1}(f)$
of $f$, which for free modules is just the ideal of $(t+1)$-minors
of a defining matrix of $f$ (see \ref{defdet} for the general
definition). We will present two types of results for these ideals
or, equivalently, the degeneracy loci
$D_t(f):=\Spec(A/I_{t+1}(f))$, namely on \bnum[(I)] \item the
dimension of $D_t(f)$, \item the connectedness of $D_t(f)$. \enum

In commutative algebra there are various results known in this
direction. Concerning (I), the classical Krull
principal ideal theorem states that the dimension drops at most
by one when we factor out an element $x\in \fm$. More generally
for determinantal ideals there is the well known Eagon-Northcott
bound on the height of such ideals \cite{EN}.

The principal ideal theorem was generalized to order ideals by
Eisenbud and Evans \cite{EE} and in its general form by Bruns
\cite{Br}. More precisely, if $N$ is a finite $A$-module, $f\in\fm
N$ and $N^\vee(f)=\{\varphi (f)| \varphi\in \Hom_A(N,A)\}$ is the
order ideal of $f$, then under suitable hypotheses the dimension
of the quotient $A/N^\vee(f)$ is at least $\dim A-\rk N$. In the
setup above this is the case where $M=A$ and $t=0$, i.e.\
$N^\vee(f)=I_1(f)$. In the meantime there has been a wide range of
extensions and applications of these results, see e.g.\
\cite{EHU1, EHU3, EHU4}.

Concerning (II), much less is known. Using methods from formal
geometry Grothendieck showed in \cite[Exp.\ XIII]{SGA} the
analogue of the principal ideal theorem for connectedness.
Grothendieck's theorem says in its simplest form that for a
noetherian complete local domain $A$ of dimension $\ge 3$ and
$x\in \fm$ the punctured spectrum $\Spec( A/xA)\backslash\{\fm\}$
is connected (see \cite{FOV}, Section 3.1 for a systematic
treatment of connectedness and a more elementary proof of this
result). In a more recent paper F.\ Steffen \cite{St1} generalized
this to determinantal ideals of matrices; however her method is
restricted to characteristic 0.

Our main results are as follows. First we show: {\em if $f\in
\fm\Hom_A(M,N)$ and if for some minimal prime $\fq$ of $A$ of
maximal dimension the modules $M_\fq$ and $N_\fq$ are free over
$A_\fq$ of ranks $m$ and $n$, respectively, then $\dim D_t(f)\ge
\dim A-\tau$ with $\tau:=(n-t)(m-t)$}. This generalizes the
Bruns-Eisenbud-Evans principal ideal theorem for order ideals. If
moreover {\em $M$ and $N$ are free at all primes $\ne\fm$ and
$\hat A$ is connected in dimension $d$, then $D_t(f)$ is connected
in dimension $d-\tau$}. Applying this to free modules we recover
the connectedness theorem of F.\ Steffen mentioned above in any
characteristic. Recall that a noetherian scheme $Y$ is said to be
$\it connected \ in \ dimension  \ d$ for $d$ an integer, if
$\dim Y>d$ and if $Y\setminus X$ is connected for every closed
subset $X\subseteq Y$ with $\dim X < d$.

In projective geometry Fulton and Lazarsfeld showed analogous
results for maps of vector bundles $f:\cM\to \cN$ if
$\cHom(\cM,\cN)$ is ample.  Ideally the local statements should
imply the corresponding projective results by passing to the
affine cones. However, the assumption that $f\in \fm\Hom(M,N)$
turns out to be too restrictive for such applications. We are able
to obtain such results in characteristic $p>0$ provided that
$\cHom(\cM,\cN)$ is cohomologically $p$-ample in the sense of
Kleiman \cite{Kl}. To deduce this from a local statement we
introduce in Section\ 2 as a new tool the concepts of
(cohomological) $p$-ampleness and ampleness for modules over local
rings that reflect to a certain extent the corresponding notions in
projective geometry, see Propositions \ref{pample} and
\ref{ample}. This allows us to weaken the assumption ``$f\in
\fm\Hom(M,N)$" to the condition that the pair $(Af, \Hom(M,N))$ is
$p$-ample, see Theorems \ref{gencodim} and \ref{gencon}.

The $A$-dual of an $A$-module $M$ will be denoted by $M^\vee$. If
$(A,\fm)$ is local we will say that $M$ is a module of {\em rank
m} with an {\em isolated singularity} if $M_\fp$ is free of
constant rank $m$ for all non-maximal primes $\fp$ of $A$. As a
general reference for basic facts on commutative algebra and
algebraic geometry we refer to the books of Eisenbud \cite{Ei},
Matsumura \cite{Ma} and Hartshorne \cite{Ha1}.

\section{$p$-Ampleness and ampleness over local rings}

In this section we propose a definition of
$p$-ampleness and ampleness over local rings.
Such notions should reflect the usual $p$-ampleness and
ampleness in global geometry as introduced by Hartshorne
\cite{Ha}. We recall shortly the relevant definitions in the projective
case.

Let $R=R_0[R_1]$ be a finitely generated standard graded algebra
over a field $K=R_0$ and $M$ a finitely generated graded $R$
module such that the associated sheaf $\cM=\tilde M$ on $X:=\Proj
R$ is locally free. Then $\cM$ is called {\em ample} if for every
other coherent sheaf $\cN$ on $X$ the sheaf $\S^n(\cM)\otimes_X
\cN$ is globally generated for $n\gg 0$, where $\S^n$ denotes the
$n$th symmetric power functor. We note that this definition is
equivalent to the cohomological condition
$$
H^j(X,\S^n(\cM)\otimes_X \cN)=0 \quad\mbox{for all }j>0,\,\, n\ge
n_0=n_0(\cN)\leqno (\sharp)
$$
and all $\cN$ as above, see \cite[3.3]{Ha}.

If $K$ has characteristic $p>0$, $p$-{\em ampleness} means by
definition that for every coherent sheaf $\cN$ on $X$ the sheaf
$F^{a\star}(\cM)\otimes_X \cN$ is globally generated for $a\gg 0$,
where $F:X\to X$ is the Frobenius map, see \cite[Sect.\ 6]{Ha}. We
note that one cannot characterize this property by a cohomological
condition similar to the one in $(\sharp)$. Thus S.\ Kleiman
\cite{Kl} introduced the notion of {\em cohomologically
$p$-ample}\/ locally free sheaves meaning that
$$
H^j(X,F^{a\star}(\cM)\otimes_A\cN)=0 \quad\mbox{for all }j>0,\,\,
n\ge n_0=n_0(\cN)
$$
and all $\cN$ as above.

The following implications are known: 
\bdi \mbox{cohomologically
$p$-ample} & \rImplies^{\mbox{\cite{Kl}}}& \mbox{$p$-ample}&
\rImplies^{\mbox{\cite[6.3]{Ha}}}& \mbox{ample.} 
\edi 
However, the converse implications do not hold as shown by Kleiman
\cite{Kl} and  Gieseker \cite[p.\ 111]{Gi}.

These notions heavily depend on the grading. Thus to define a
local analogue our idea is to consider a relative notion depending
on a pair $(M',M)$, where $M'$ is a submodule of $M$.

\bsit 
Let $(A,\fm)$ be a noetherian local ring of prime
characteristic $p>0$ and $M$ an $A$-module. The Frobenius map
$F:A\to A':=A$ is a ring homomorphism, and
$$
M\mapsto F(M):=M\otimes_F A=M\otimes_{A}A'
$$
is a functor on the category of $A$-modules. Iterating $a$ times
we obtain the map $F^a:A\to A':=A$ and the corresponding functor
$F^a$. In order to distinguish between duals into $A$ and $A'$ we
write $M^\vee=\Hom_{A}(M,A)$ and $M^\times=\Hom_{A'}(M,A')$.
\esit

For our purposes the following definition of $p$-ampleness is the
most appropriate:

\bdfn\label{2.2}
A pair $(M',M)$, where $M'\subseteq M$ are
$A$-modules, is called {\em $p$-ample} if for $a\gg 0$
the image of the natural map
$$
F^a(M')\to F^a(M)^{\times\times}
$$
is contained in $\fm F^a(M)^{\times\times}$. Similarly, for
$R=R_0[R_1]$ a noetherian graded algebra over a local ring $A=R_0$
and $M' \subseteq M$ graded $R$-modules we say that the pair
$(M',M)$ is $p$-ample, if the condition above is satisfied over $R$ with $\fm$
the maximal homogeneous ideal of $R$ or, equivalently, if the pair
$(M'_\fm,M_\fm)$ is $p$-ample over the localized ring $A=R_\fm$.
\edfn

We first record some basic properties of this notion:

 \brem\label{rempample}
 \bnum[(1)] \item Any pair $(M',M)$ with $M'\subseteq\fm M$ is
$p$-ample. \item If the pair $(M',M)$ is $p$-ample, then so is
$(N,M)$ for every submodule $N$ of $M'$. \item A pair $(M',M)$ of
finite $A$-modules is $p$-ample if and only if the pair of $\hat
A$-modules $(\hat M',\hat M)$ is $p$-ample. \item Let $f:M_1\to
M_2$ be a homomorphism of $A$-modules. If $(M',M_1)$ is $p$-ample,
then so is $(f(M'),M_2)$. \item If $(M',M)$ is $p$-ample, then
$(N\otimes_A M', N\otimes_A M)$ is $p$-ample for any $A$-module
$N$. \item If $(M',M)$ is $p$-ample, then so is $(M'\cdot
\S^{n-1}M, \S^nM)$ for every $n> 0$. \enum \erem

\bproof Parts (1), (2), (3) and (4) are immediate from the
definitions. The assertions (5) and (6) are an easy consequence of
the fact that the Frobenius functor is compatible with tensor
products and symmetric powers. \eproof

In the following proposition we describe the relationship of our
notion of $p$-ample\-ness with the corresponding projective
notion.

\bpro\label{pample} Let $R=R_0[R_1]$ be a finitely generated
graded algebra over a field $K=R_0$ of characteristic $p>0$ with
$\depth R\ge 2$ and let $M=\bigoplus_{i\in \Z}M_i$ be a finite
graded $R$-module with an isolated singularity. Then the following
hold: \bnum[$(1)$] \item If the associated locally free sheaf
$\cM:=\tilde M$ on $X:=\Proj R$ is cohomologically $p$-ample, then
the pair $(R\cdot M_0 ,M)$ is $p$-ample. \item If $\cM$ is
generated by the sections in $M_0$ and if the pair $(R\cdot M_0
,M)$ is $p$-ample, then $\cM$ is $p$-ample on $X$. \enum \epro

\bproof
Since $\depth A\ge 2$ and $M$ has an isolated singularity we
have
$$
F^a(M)^{\times\times}\cong
H^0_\bullet(F^{a\star}(\cM)):=\bigoplus_{i\in \Z}H^0(X,
F^{a\star}(\cM)(i)).\leqno(*)
$$
In order to show (1) let $x_1,\ldots, x_k\in R_1$ be a $K$-basis.
The Koszul complex $K(\ul{x}, F^{a}(M))$  on $\ul{x}$ provides an
exact sequence
$$
0\longrightarrow F^{a\star}(\cM)(-k)\longrightarrow \cdots
\longrightarrow F^{a\star}(\cM)(-1)^{\binom{k}{1}} \longrightarrow
F^{a\star}(\cM)\to 0\,
$$
on $X$. As $\cM$ is cohomologically $p$-ample the groups
$H^j(X,F^{a\star}(\cM)(-i))$ vanish for $j> 0$ and $0\le i\le k$,
if $a\gg 0$. Thus taking sections the sequence
$$
  \cdots \longrightarrow H^0(X, F^{a\star}(\cM)(-1))^{\binom{k}{1}}
\longrightarrow H^0(X,F^{a\star}(\cM))\longrightarrow 0
$$
is exact. In view of $(*)$ this implies that $  R_1\cdot
F^{a}(M)^{\times\times}_{-1}= F^{a}(M)^{\times\times}_0
\supseteq F^{a}(RM_0)_0$ as desired.

If the conditions in (2) are satisfied then
$$
F^{a}(RM_0)_0\subseteq R_1F^{a}(M)^{\times\times}_{-1}=R_1 H^0(X,
F^{a\star}(\cM)(-1))
$$
for $a\gg 0$. As $\cM$ is globally generated by $M_0$ it follows
that the map
$$
\mathcal{O}_X(1)\otimes H^0(X,F^{a\star}(\cM)(-1)) \longrightarrow
F^{a\star}(\cM)
$$
is surjective. Hence $F^{a\star}(\cM)$ is a quotient of a direct
sum of copies of $\cO_X(1)$, and $\cM$ is $p$-ample as required.
\eproof

Although not needed in the sequel it is interesting to note that
one can also introduce a local notion of ampleness. For the
remainder of this section let $(A,\fm)$ be an arbitrary noetherian
local ring.

\bdfn A pair $(M',M)$, where $M'\subseteq M$ are $A$-modules, is
called {\em ample} if for $n\gg 0$ the image of the natural map
$\S^n M'\to(\S^nM)^{\vee\vee}$ is contained in $\fm
(\S^nM)^{\vee\vee}.$ \edfn

As in Definition \ref{2.2} one can define a similar notion in the
graded case. In analogy to \ref{rempample} we note some basic
properties of this notion.

\brem\label{remample} 
\bnum[(1)] 
\item Any pair $(M',M)$ with
$M'\subseteq\fm M$ is ample. 
\item If the pair $(M',M)$ is ample,
then so is $(N,M)$ for every submodule $N$ of $M'$. 
\item A pair
$(M',M)$ of finite $A$-modules is ample if and only if the pair of
$\hat A$-modules $(\hat M',\hat M)$ is ample. 
\item Let $f:M_1\to
M_2$ be a homomorphism of $A$-modules. If $(M',M_1)$ is ample then
so is $(f(M'),M_2)$. 
\item If $\Q\subseteq A$ and $(M',M)$ is
ample, then $(M'\otimes_A N, M\otimes_A N)$ is ample for any
$A$-module $N$. 
\item
If $\Q\subseteq A$ and $(M',M)$ is ample, then so is $(M'\cdot
\S^{n-1}M, \S^nM)$ for every $n> 0$. 
\enum \erem

\bproof The proofs of (1), (2), (3) and (4) are immediate from the
definitions. Part (5) follows from the fact that in characteristic
0 there is a natural surjective map $N^{\otimes n}\otimes_A
\S^n(M)\to
\S^n(N\otimes_A M)$. Finally, part (6) is an easy consequence of
(5) and (4), as the pair $(M'\otimes_A
\S^{n-1}M, M\otimes_A\S^{n-1}M)$ is ample and maps onto the pair
$(M'\cdot \S^{n-1}M, \S^nM)$. We leave the simple details to the
reader. \eproof

\brem By a well-known result of Barton \cite{Ba} the tensor
product of ample vector bundles on a projective variety is again
ample in any characteristic. Therefore one should expect that
\ref{remample}(5),(6) remain valid in positive characteristic as
well. \erem

The following proposition explains our terminology of ample
pairs. Its proof follows along the same line of arguments as in
\ref{pample} and is left to the reader.

\bpro\label{ample} Let $R=R_0[R_1]$ be a finitely generated graded
algebra over a field $K=R_0$ with $\depth R\ge 2$ and let
$M=\bigoplus_{i\in \Z}M_i$ be a finite graded $R$-module with an
isolated singularity. Then the following hold $:$ \bnum[$(1)$]
\item If the associated locally free sheaf $\cM:=\tilde M$ on
$X:=\Proj R$ is ample, then the pair $(R\cdot M_0 ,M)$ is ample.
\item If $\cM$ is generated by the sections in $M_0$ and if the
pair $(R\cdot M_0 ,M)$ is ample then $\cM$ is ample. \enum \epro

Let $A$ be a noetherian local ring, let $M$ be a finitely
generated $A$-module, and consider an embedding of $M^{\vee\vee}$
into a free $A$-module $F$. We recall that the {\em integral
closure} of a submodule $N\subseteq M$ is the set of all elements
in $M$ whose image in $F$ is integral over the image of the
symmetric algebra $\S(N)$ in $\S(F)$ (this definition does not
depend on the chosen embedding $M^{\vee\vee}\subseteq F$, see
\cite{EHU2}). We have the following relation with our notion of
ampleness:

\bpro Let $(A,\fm)$ be a noetherian local ring and let
$M'\subseteq M$ be finite $A$-modules. Assume that $M_\fq$ is free
for all associated primes $\fq$ of $A$. If $M'\subseteq\ol{\fm M}$
is in the integral closure of $\fm M$ in $M$, then the pair
$(M',M)$ is ample \epro

\bproof Write $M^n$ for the image of $\S^nM$ in $\S^nF$. As $M$ is
free at the associated primes of $A$ we have
$(\S^nM)^{\vee\vee}=(M^n)^{\vee\vee}$. By our assumption, the
image of $\S^nM'$ in $M^n$ is contained in $\fm M^n$ for $n \gg
0$, and hence is contained in $\fm
(M^n)^{\vee\vee}=\fm(\S^nM)^{\vee\vee}$. \eproof

\section{Main technical results}

\begin{sit}\label{1.0} Let $A$ be a ring, $E$ an $A$-module and
$f\in E^\vee$. The linear form $f:E\to A$ extends uniquely to a
morphism of $A$-algebras denoted
$$
\varphi_f:\S(E)\lto A,
$$
where $\S(E)=\bigoplus_{i=0}^\infty (\S^iE)T^i$ is the symmetric
algebra of $E$ over $A$.  Note that if $E_\fp$ is free for some
prime $\fp\in \Spec A$, then the localization of the symmetric
algebra $\S_{A_\fp}(E_\fp)$ is isomorphic to a polynomial ring
over $A_\fp$ in $\rk (E_\fp)$ indeterminates.
\end{sit}

The next theorem is one of our main technical results:

\bthm\label{1.1} Let $(A,\fm)$ be a noetherian local ring, $E$ a
finite $A$-module and $f\in E^\vee$. Assume that
$$
h \circ f\in \fm'\Hom_A(E,A') \leqno(*)
$$
for some local homomorphism $h: A \rightarrow A'$ with $(A',\fm')$
a noetherian local ring and $\dim A'/\fm A' =0$. Further suppose
that there exists a prime ideal $\fq'$ of $A'$ with contraction
$\fq:=h^{-1}(\fq')$ so that $\dim A'/\fq'=\dim A$ and $E_{\fq}$ is
a free $A_{\fq}$-module of rank $r$. Let $I\subseteq S:=\S(E)$ be
a homogeneous ideal with $I\cap A \subseteq \fq$ and $\dim
S_\fq/I_\fq=r-\tau$. If $I(f):=\varphi_f(I)$
then $\dim A/I(f) \ge \dim A-\tau$.
\ethm

We will reduce Theorem \ref{1.1} to the special case where $h$ is
the identity map:

\bthm\label{1.2} Let $(A,\fm)$ be a noetherian local ring, $E$ a
finite $A$-module and assume that
$$
f\in \fm E^{\vee}. \leqno(**)
$$
Further suppose there exists a prime ideal $\fq$ of $A$ so that
$\dim A/\fq=\dim A$ and $E_{\fq}$ is a free $A_{\fq}$-module of
rank $r$. Let $I\subseteq S:=\S(E)$ be a homogeneous ideal with
$I\cap A \subseteq q$ and $\dim S_\fq/I_\fq=r-\tau$. If
$I(f):=\varphi_f(I)$
then $\dim A/I(f) \ge \dim A-\tau$.
\ethm

\begin{proof}[Proof of Theorem $\ref{1.1}$]
We show how Theorem \ref{1.1} can be deduced from Theorem
\ref{1.2}. For this we consider $E':=E\otimes_AA'$, $f':=f\otimes
1 \in E'^{\vee}$, $S':=S\otimes_AA'=\S(E')$, and the extended
ideal $I':=IS'$. Notice that $I'(f')=I(f)A'$. If $t_1,\ldots,
t_n\in\fm$ is a system of parameters of $A$, then $\dim
A'/(t_1,\ldots,t_n)A'=\dim A'/\fm A'=0$ and so by Krull's
principal ideal theorem $\dim A'\le n=\dim A$. Likewise we have
$\dim A'/I'(f')
\leq \dim A/I(f)$. Since $\dim A'/\fq' \leq \dim
A' \leq \dim A$ and the outer terms are equal, it follows that
$q'$ is a minimal prime of $A'$ having maximal dimension and that
$\dim A'=\dim A$. As $S'_{\fq'}$ and $S_{\fq}$ are polynomial
rings in the same $r$ variables over the Artinian local rings
$A'_{\fq'}$ and $A_{\fq}$, reducing modulo the nilradicals of the
latter rings one sees that $\dim S'_{\fq'}/I'_{\fq'}=\dim
S_{\fq}/I_{\fq}=r-\tau$. Thus Theorem \ref{1.2} gives $\dim
A'/I'(f')\geq \dim A' - \tau$. Therefore $$ \dim A/I(f) \geq \dim
A'/I'(f') \geq \dim A' - \tau = \dim A - \tau,$$ 
proving Theorem \ref{1.1}.\end{proof}

It remains to prove Theorem \ref{1.2}. For this we will use the
following setting:

\bsit\label{1.2a} Fix a set of generators $f_1,\ldots,f_k$ for
$E^\vee$ and let
$$
\varphi:S=\S(E)\to A[\ul{Y}]:=A[Y_1,\ldots, Y_k]
$$
denote the (unique) map of $A$-algebras with
$$
\varphi(eT)=f(e)+\sum_{i=1}^k f_{i}(e)Y_i
$$
for $e\in E$. Note that $\varphi$ is not homogeneous in general
when we equip $\S(E)$ and $A[\ul{Y}]$ with their standard grading.
In the following let $J:=IA[\ul{Y}]$ be the extension ideal, $\fM$
the maximal homogeneous ideal $\fm A[\ul{Y}]+(\ul{Y})$ of
$A[\ul{Y}]$ and $B:=A[\ul{Y}]_\fM$ the localization. Obviously
$\varphi_f$ is the composed map
\begin{diagram}
S  & \rTo^{\varphi} & A[\ul{Y}] &
\rTo^{Y_i=0} & A\,.
\end{diagram}
\esit

The next lemma is a standard fact; for completeness we give a
short argument.

\begin{lem}\label{1.3} If for some prime $\fp$ of $A$ the
localized module $E_\fp$ is free of rank $r$, then the maps \bdi
A_\fp&\rTo^{can}& S_\fp&\rTo^\varphi& A_\fp[\ul{Y}] \edi are
polynomial ring extensions of relative dimension $r$ and $k-r$,
respectively.
\end{lem}

\begin{proof} Localizing at $\fp$ we may suppose
that $E$ is free over $A$ and so $E^\vee$ has a basis, say
$f_1,\ldots,f_r$, where $r=\rk E$. We can write
$f=\sum_{\varrho=1}^r a_{\varrho} f_\varrho$ and
$f_i=\sum_{\varrho=1}^r a_{i\varrho} f_\varrho$ for $i >r$ with
$a_\varrho,a_{i\varrho}\in A$. The basis $f_1,\ldots,f_r$ provides
an isomorphism between the symmetric algebra $S$ and the
polynomial ring $A[Y_1,\ldots,Y_r]$, and the map $S\cong
A[Y_1,\ldots, Y_r] \to A[\ul{Y}]$ is given by $Y_\varrho\mapsto
a_\varrho +Y_\varrho+ \sum_{i>r} a_{i\varrho} Y_i$,\ \ $1\le i\le
r$. Hence via this map, $A[\ul{Y}]$ may be considered as a
polynomial ring over $S$ in the indeterminates $Y_{r+1},\ldots,
Y_k$.
\end{proof}

The following lemma is a crucial step in the proof of Theorems
\ref{1.2} and \ref{1.6}; it is the only place where we use the
assumption $(**)$.

\blem\label{1.4} With the notation of $\ref{1.2a}$ and the
assumptions of Theorem $\ref{1.2}$ the following hold $:$

\bnum[$(1)$] \item All minimal primes of $J\subseteq A[\ul{Y}]$
are contained in $\fM$. \item Assume moreover that $A=\hat A$ is
complete and $E$ has an isolated singularity. If $S/I:\fm^\infty$
is an integral domain then so is the completion of
$B/JB:\fm^\infty$. \enum \elem

\bproof As $f\in \fm E^\vee$ we can write $f=\sum_it_if_i$ with
$t_i\in \fm$. Let $\psi:A[\ul{Y}]\to A[\ul{Y}]$ be the isomorphism
of $A$-algebras given by $Y_i\mapsto Y_i-t_i$. The composition
$$
\psi\circ\varphi:S\longrightarrow A[\ul{Y}]
$$
maps $eT$ onto $f(e)+\sum f_i(e)(Y_i-t_i)=\sum f_i(e)Y_i$ and so
$\psi\circ\varphi$ is homogeneous. The homogeneous ideal
$\psi\circ\varphi(I)A[\ul{Y}]=\psi(J)$ has only homogeneous
minimal primes that are all contained in $\fM$. As $\psi$ maps
$\fM$ into itself, (1) follows.

In order to prove (2) we first notice that for every prime
$\fp\ne\fm$ of $A$ the map $\psi\circ\varphi: S_\fp/I_\fp\to
A_\fp[\ul{Y}]/\psi(J_\fp)$ is a polynomial ring extension by Lemma
\ref{1.3}. Thus it follows that $A[\ul{Y}]/(\psi(J):\fm^\infty)$
is a homogeneous domain and so its completion, which is equal to
the completion of $B/(\psi(J)B:\fm^\infty)$, is a domain. As the
rings $B/(JB:\fm^\infty)$ and $B/(\psi(J)B:\fm^\infty)$ are
isomorphic under the map induced by $\psi$, (2) follows. \eproof

\begin{proof}[Proof of Theorem $\ref{1.2}$]
We may assume that $A=\hat A$ is complete. Furthermore, if $\dim
(A/(\fq+I(f))\ge \dim A/\fq-\tau$, then $\dim A/I(f)\ge \dim A-
\tau$. Thus, replacing $A$ by $A/\fq$ we are reduced to the case
where $A$ is a complete domain and $\fq=0$. Let $\fp$
be a minimal prime of $I$ that extends to a minimal prime of
$I_\fq$ of maximal dimension in $S_\fq$. Clearly it suffices to
show the result for $\fp$ instead of $I$. Thus we may further
suppose that $I$ is prime.

According to Lemma \ref{1.3}, $S_\fq/I_\fq\to A_\fq[\ul{Y}]/J_\fq$
is a polynomial extension of relative dimension $k-r$. Thus there
is a minimal prime $\fP$ of $J$ contracting to $I$. Lemma
\ref{1.4}(1) shows that $\fP$ is contained in $\fM$. Applying the
dimension formula \cite[5.5.8.1]{EGA} for morphisms of finite type
we get
$$
\dim B/\fP B=\dim A+\trdeg_{A_\fq} B_\fq/\fP B_\fq.
$$
By Lemma \ref{1.3}, $\trdeg_{A_\fq} B_\fq/\fP B_\fq=\trdeg_{A_\fq}
S_\fq/I_\fq +(k-r)=k-\tau$, where we used our assumption that
$\dim S_\fq/I_\fq=r-\tau$. Thus we obtain
$$
\dim B/JB\ge \dim B/\fP B=\dim A+k-\tau.
$$
Using Krull's principal ideal theorem it follows that
$$
B/(JB+(Y_1,\ldots,Y_k))\cong  A/I(f)
$$
has dimension $\ge \dim A-\tau $, as desired.
\end{proof}

\bexa \label{1.5a} Let $(A, \fm)$ be the localization of the ring
$R:=\C[X_0,\ldots, X_3]$ at the homogeneous maximal ideal.
Consider the $A$-module $E:=\coker (h: A\to \Theta_A)$, where
$\Theta_A$ is the module of derivations on $A$ and $h$ is given by
$h(1)=\sum X_i\partial /\partial X_i$. The element
$f:=X_0dX_1-X_1dX_0+X_2dX_3-X_3dX_2\in E^\vee$ has then the
property $f(E)=\fm $. Taking $I\subseteq \S(E)$ to be the ideal
generated by the forms of strictly positive degree we have $\dim
\S(E_{\fq})/I_\fq=0=\rk E-3$ with $\fq:=0$ the minimal prime of
$A$, but $\dim A/I(f)=0<\dim A-3$. This example shows that the
assumption $f\in\fm E^\vee$ is essential.

It is instructive to understand geometrically why Lemma
\ref{1.4}(1) does not hold in this example. The ring $A$ can be
regarded as the local ring of the affine 4-space $V:= \Spec R$ at
the origin $0\in V$ and $E$ is the localization of $E_R:=\coker
(h: R\to \Theta_R)$ with $h$  as above. Thus $\V(E_R):=\Spec
\S(E_R)$ embeds into the  affine space $V^\vee\times
V=\V(\Theta_R)$. Taking as generators of $E^\vee_R$ the standard
syzygies $f_{ij}=X_idX_j-X_jdX_i$ we have $\Spec R[\ul{Y}]\cong
\Lambda^2 V^\vee\times V:=\V(\Lambda^2 \Theta_R)$, and the
generators $f_{ij}$ induce a morphism $H:\Lambda^2 V^\vee\times
V\to\V(E_R).$ The composition with the embedding $ \V(E_R)$ into
$V^\vee\times V$ on the $\C$-valued points is the map
$$
\ba{c} \Lambda^2 V^\vee\times V\lto \V(E_R)\lto V^\vee\times
V\\
(T, x)\mapsto (Tx,x).
\ea
$$
Here we consider a $\C$-valued point of  $\Lambda^2 V^\vee$ as a
skew symmetric matrix $T$. The given element $f$ corresponds to a
nondegenerate matrix $T_0$ and the map $\varphi$ in \ref{1.2a} to
$\Phi:(T,x)\mapsto ((T+T_0)x,x).$ The ideal $I$ defines the zero
section $\{0\}\times V$, and the preimage of this zero section in
$\Lambda^2 V^\vee\times V$ is the set
$$
Z:=\{(T,x): (T+T_0)x=0\}.
$$
The irreducible components of $Z$ are given by
$$
Z_1:=\{(T,x): (T+T_0)x=0\mbox{ and }\Pf(T+T_0)=0\}
\quad\mbox{and}\quad
Z_2:=V^\vee\times \{0\},
$$
where $Z_1$ is of codimension 3 and $Z_2$ of
codimension 4. Thus $Z_1$
does not meet $\{0\}\times V$, and the hypersurfaces
$Y_{ij}=0$ intersect $Z$ in the origin
$(0,0)$.
\eexa

\begin{rem}\label{r1.7} In the setting of Theorem \ref{1.1} there is another condition
$(*$$**)$ that implies $(*)$, but is weaker than assumption $(**)$
in Theorem \ref{1.2}. We abbreviate as before
$-^\times=\Hom_{A'}(-,A')$. If
$$
f\otimes 1\in \fm'(E^{\vee}\otimes_AA')^{\times\times},
\leqno(\mbox{$**$$*$})
$$
then $(*)$ holds. This implication follows from the existence of
the natural $A'$-linear map

$$
(E^\vee\otimes_AA')^{\times\times}\longrightarrow
(E\otimes_AA')^{\times\times\times} \longrightarrow
(E\otimes_AA')^{\times} \cong\Hom_A(E,A'),
$$
that sends the image of $f\otimes 1$ to $h\circ f$.
\end{rem}

Condition $(*$$**)$ for $h$ any high power of the Frobenius map is
precisely the $p$-ampleness of the pair $(Af, E^\vee)$ as
introduced in Definition \ref{2.2}. Thus we deduce:

\bthm\label{1.2n} Assume that $(A,\fm)$ is a noetherian local ring
of prime characteristic $p$. The conclusion of Theorem $\ref{1.2}$
remain valid if one replaces assumption $(**)$ by the weaker
condition that the pair
$$
(Af, E^\vee)\quad\mbox{is $p$-ample}. \leqno (\mbox{$**$$**$})
$$
\ethm

We now turn to our results about connectedness.

\bthm\label{1.6} Let $(A,\fm)$ be a noetherian local ring, $E$ a
finite $A$-module of rank $r$ with an isolated singularity, and
$f\in E^\vee$. Assume that
$$
h \circ f\in \fm'\Hom_A(E,A') \leqno(*)
$$
for some local homomorphism $h: A \rightarrow A'$ with $(A',\fm')$
a noetherian local ring so that the map of completions $\hat
h:\hat A \rightarrow \hat A'$ has nilpotent kernel and $\hat{A'}$
is integral over the image of $\hat{h}$.
Let $I\subseteq S:=\S(E)$ be a homogeneous ideal and suppose that
the map
$$
\Spec (S/I)\backslash V(\fm(S/I))\lto \Spec (A)\backslash \{\fm\}
$$
has non empty geometrically irreducible fibers all of the same
dimension $r-\tau$. If $\hat A$ is connected in dimension $d$ then
$\hat A/I(f)\hat A$ is connected in dimension $d-\tau $, where
$I(f):=\varphi_f(I)$. \ethm

Once again, the special case $h={\rm id}_A$ is noteworthy:

\bthm\label{1.6a} Let $(A,\fm)$ be a noetherian local ring, $E$ a
finite $A$-module of rank $r$ with an isolated singularity, and
assume that
$$
f\in \fm E^{\vee}. \leqno(**)
$$
Let $I\subseteq S:=\S(E)$ be a homogeneous ideal and suppose that
the map
$$
\Spec (S/I)\backslash V(\fm(S/I))\lto \Spec (A)\backslash \{\fm\}
$$
has non empty geometrically irreducible fibers all of the same
dimension $r-\tau$. If $\hat A$ is connected in dimension $d$ then
$\hat A/I(f)\hat A$ is connected in dimension $d-\tau $, where
$I(f):=\varphi_f(I)$. \ethm

In the light of Remark \ref{r1.7}, applying Theorem \ref{1.6} to
the iterated Frobenius $h=F^a$,
$a\gg 0$, gives:

\bthm\label{1.6n} Assume that $(A,\fm)$ is a noetherian local ring
of prime characteristic $p$. The conclusion of Theorem $\ref{1.6}$
remains valid if one replaces assumption $(**)$ by the weaker
condition that the pair
$$
(Af, E^\vee)\quad\mbox{is $p$-ample}. \leqno (\mbox{$**$$**$})
$$
\ethm

To prove Theorem \ref{1.6} we will need the next two lemmas.

\blem\label{lemcon} Let $Y$ be a noetherian scheme and let
$Y_i\subseteq Y$ be closed subsets of dimension $\ge d+1$ with
$Y=Y_1\cup \ldots\cup Y_s$. Consider the following conditions $:$
\bnum[$(i)$] \item $Y$ is connected in dimension $d$. \item For
all $Y_i$, $Y_j$ there exists a sequence $Y_{i_1}, \ldots,
Y_{i_{t}}$ such that $Y_i=Y_{i_1}$, $Y_j=Y_{i_t}$ and
$$
\dim(Y_{i_\alpha }\cap Y_{i_{\alpha +1}}) \ge d\,\, \mbox{ for
}\,\, 1\le \alpha  < t.
$$
\enum Then $(i)$ implies $(ii)$. If all $Y_i$ are connected in
dimension $d$ then also the converse is true. \elem

\bproof For (ii)$\Rightarrow$(i) we refer the reader to
\cite[Proposition 3.1.4]{FOV}. To prove (i)$\Rightarrow$(ii), note
first that the relation given by $i\sim j$ if and only if (ii) is
satisfied, is an equivalence relation. For an equivalence class
$\alpha\subseteq \{1,\ldots,s\}$ let $Y_\alpha:=\bigcup_{i\in
\alpha} Y_i$. Clearly $Y=\bigcup_\alpha Y_\alpha$ and $\dim
(Y_\alpha\cap Y_\beta )< d$ for $\alpha\ne\beta$. Thus $Y$ will be
not connected in dimension $d$ unless there is only one
equivalence class. \eproof

\blem\label{1.7} Let $Y$ be an irreducible noetherian scheme and
$\cR$ a finitely generated  positively graded $\cO_Y$-algebra with
coherent homogeneous components. If the morphism $ X:=\Spec \cR
\lto Y$ has irreducible fibers all of the same dimension, then $X$
is irreducible. \elem

\bproof If $\cR$ is a finite $\cO_Y$-module then the assertion is
immediate. Otherwise the morphism
$$
\pi: Z:=\Proj \cR\longrightarrow Y
$$
has non empty irreducible fibers all of the same dimension, say
$s\ge 0$. Let $Z_0$ be an irreducible component of $Z$ with
generic point $\eta_0$ such that $\pi(Z_0)=Y$. The generic fiber
$F:= \pi^{-1}(\pi(\eta_0) )$ through $\eta_0$ contains $F\cap Z_0$
as an irreducible component and so these sets are equal by our
assumption. By the dimension formula $\dim Z_0=s+\dim Y$. Using
the semicontinuity of fiber dimensions all fibers of $Z_0\to Y$
have dimension $\ge s$. Therefore they are equal to the fibers of
$\pi$ as these are irreducible of dimension $s$ by our assumption.
Hence $Z=Z_0$ and so $Z$ and then also $X$ are irreducible, as
desired. \eproof

\bproof[Proof of Theorem $\ref{1.6}$.] We may assume that $A=\hat
A$ and $A'=\hat A'$ are complete. We may further suppose that
$d>0$. Let $Y_{i}$ be the irreducible components of $Y=\Spec A$
and suppose that we can prove that $V_{i}:=Y_i\cap V$ is connected
in dimension $d-\tau$, where $V:=V(I(f))$. We claim that then $V$
is connected in dimension $d-\tau$. Indeed, applying Lemma
\ref{lemcon}, for all $Y_i$, $Y_j$ there exists a sequence of
irreducible components $Y_i= Y_{i_1}, \ldots, Y_{i_{t}}= Y_j$,
with $\dim(Y_{i_\alpha }\cap Y_{i_{\alpha +1}}) \ge d$ for $1\le
\alpha < t$. Taking the intersection with $V$ and using Theorem
\ref{1.1} we obtain that $\dim(V_{i_\alpha }\cap V_{i_{\alpha
+1}}) \ge d-\tau$ for $1\le \alpha  < t$. Hence again by Lemma
\ref{lemcon}, $V$ is connected in dimension $d-\tau$.

As $h: A \to A'$ has nilpotent kernel and $A'$ is integral over
the image of $h$, the induced map $\Spec A' \to \Spec A$ is
surjective. Thus, by factoring out an arbitrary minimal prime of
$A$ and a prime of $A'$ lying over it, we may suppose that $A
\subseteq A'$ is an integral extension of domains. As $\dim
A'=\dim A
> d$ it follows that $A'$ is connected in dimension $d$.
With $I'(f')$ defined as in the proof of Theorem \ref{1.1} one has
$I'(f')=I(f)A'$. Since $A \subseteq A'$ is an integral extension,
connectedness in dimension $d$ descends from $A'/I(f)A'$ to
$A/I(f)$. Replacing $A$ by $A'$ we may then assume that $A=A'$ and
that the stronger condition $(**)$ of Theorem \ref{1.6a} holds.
Furthermore $A$ is a complete domain.

After replacing the ideal $I$ by its radical, Lemma \ref{1.7}
shows that $S/I:\fm^\infty$ is a domain, and then by Lemma
\ref{1.4}(2) the same is true for the completion of
$B/JB:\fm^\infty$. Applying Grothendieck's connectedness theorem
we obtain that
$$
V(JB:\fm^\infty+(Y_1,\ldots, Y_k))\subseteq \Spec B/(Y_1,\ldots,
Y_k)\cong \Spec A
$$
is connected in dimension $d-\tau $. As this set is equal to
$V(I(f))$ the result follows.
\eproof

We can also generalize our results in the spirit of \cite[Theorem
3.1]{EHU3}. For instance:

\bcor\label{1.8}  Theorem $\ref{1.2}$ remains valid if we replace
condition $(**)$ by the weaker condition that \bnum \item[$(**')$]
there exists a local map $h:A\to A'$ to a noetherian local ring
$(A',\fm')$ so that $\fq=h^{-1}(\fq')$ for some prime $\fq'$ of
$A'$ and $f\otimes 1\in \fm'(E^\vee\otimes_AA'/{\rm torsion})$.
\enum \ecor

\bproof Inspecting the proof of Theorem \ref{1.2} we only need to
show the following analogue of Lemma \ref{1.4}(1):  All minimal
primes $\fP$ of $J\subseteq A[\ul{Y}]$ with $\fP\cap A=\fq$ are
contained in $\fM$.

 As $f\otimes 1\in \fm' (E^\vee\otimes_AA'/{\rm torsion})$
we can write $f\otimes 1=\sum_it_if_i\otimes 1$ with $t_i\in
\fm'$. Let $\psi:A[\ul{Y}]\to A'[\ul{Y}]$ be the morphism given by
$Y_i\mapsto Y_i-t_i$. As in the proof of Lemma \ref{1.4}(1) the
composition  $\psi\circ\varphi:S\to A'[\ul{Y}]$ is homogeneous. The
homogeneous ideal
$\psi\circ\varphi(I)A'[\ul{Y}]=\psi(J)A'[\ul{Y}]$ has only
homogeneous minimal primes that are contained in the maximal
homogeneous ideal $\fM'$ of $A'[\ul{Y}]$. Moreover $\psi$ maps
$\fM$ into $\fM'$. Since $A_{\fq}$ is Artinian and the map
$A_{\fq} \to A'_{\fq'}$ is local, it follows that every minimal
prime of $JA_{\fq}[\ul{Y}] \subseteq A_{\fq}[\ul{Y}]$ is
contracted from a minimal prime of $\psi(J)A'_{\fq'}[\ul{Y}]
\subseteq A'_{\fq'}[\ul{Y}]$. Thus for every minimal prime $\fP$
of $J\subseteq A[\ul{Y}]$ with $\fP\cap A=\fq$ there exists a
minimal prime $\fP'$ of $\psi(J) A'[\ul{Y}]\subseteq A'[\ul{Y}]$
lying over it. As $\fP' \subseteq \fM'$ we conclude that
$\fP\subseteq \fM$.
\eproof

\section{Codimension and connectedness of degeneracy loci}

\subsection{Maps of modules}
While it is clear how to define determinantal ideals and
degeneracy loci for morphisms of locally free sheaves, for maps of
arbitrary coherent sheaves several definitions are possible. For
our purposes the following natural notion is best suited.

\bsit\label{defdet} Given a morphism $f:\cM\to \cN$ of coherent
sheaves on a scheme $X$ we can introduce the $(t+1)$-$st$ {\em
determinantal ideal} of $f$ by
$$
\cI_{t+1}(f):=\mbox{Im}(\Lambda^{t+1}\cM\otimes_{\cO_{X}}
\Lambda^{t+1}\cN^\vee\lto \cO_{X}),
$$
where the map is defined by
$$
m_{0}\otimes\ldots\otimes m_{t}\otimes
\varphi_{0}\otimes\ldots\otimes \varphi_{t}\longmapsto
\det(\varphi_{i}\circ f(m_{j}))
$$
for local sections $m_{0},\ldots,m_{t}$ in $\cM$ and
$\varphi_{0},\ldots,\varphi_{t}$ in $\cN^\vee$. Its set of zeros
$D_t(f):=V(I_{t+1}(f))$ will be called the {\em $t$-$th$
degeneracy locus} of $f$.

Similarly, if $f:M\to N$ is a linear map of $A$-modules then
$I_{t+1}(f)$ denotes the image of the induced map $\Lambda^{t+1}M
\otimes_A \Lambda^{t+1}N^\vee\lto A$ and
$D_t(f):=V(I_{t+1}(f))\subseteq \Spec A$ the $t$-th degeneracy
locus.  Note that the sheaf $\tilde I_{t+1}(f)$ associated to
$I_{t+1}(f)$ on $X=\Spec A$ is just the $(t+1)$-st determinantal
ideal associated to $\tilde f:\tilde M\to\tilde N$. \esit

If $\cM$ and $\cN$ are free $\cO_{X}$-modules of rank $m$
and $n$, respectively, then $f$ can be described by a $n\times
m$-matrix $F$. In this case $I_{t+1}(f)=I_{t+1}(F)$ is just the
ideal generated by the $(t+1)$-minors of $F$.

\bthm \label{gencodim} Let $(A,\fm)$ be a noetherian local ring
and $f:M\to N$ a linear map of finite $A$-modules. Assume there
exists a prime ideal $\fq \subseteq A$ with $\dim A=\dim A/\fq$
such that $M_\fq$ and $N_\fq$ are free $A_\fq$-modules of ranks
$m$ and $n$, respectively. If $t\le\min\{n,m\}$ then with
$\tau:=(n-t)(m-t)$ we have $\dim D_t(f)\ge \dim A-\tau$ in each of
the following cases $:$ \bnum[$(a)$] \item  $f\in \fm \Hom(M,N)$.
\item $A$ has prime characteristic $p$ and the pair
$(Af,\Hom(M,N))$ is p-ample. \enum \ethm

\bproof With $E:=\Hom(M,N)^\vee$ the morphism $f$ can be
considered as an element in $E^\vee$. The idea is to apply
Theorems \ref{1.2} and \ref{1.2n} to the ``generic" determinantal
ideal $I_{t+1}\subseteq S:=\S(E)$.

In the case that $N$ is free the definition of $I_{t+1}$ is easy:
the evaluation map $\Hom(M, N)\otimes_A M\to N$ induces a map
$M\to \Hom(M, N)^\vee \otimes_A N =E\otimes_A N$. It gives a
``universal" morphism
$$
F: \S(E)\otimes_A M\longrightarrow \S(E)\otimes_A N,
$$
and one can take $I_{t+1}:=I_{t+1}(F)$. If $N$ is not free no such
universal map exists, but one can nevertheless define the generic
determinantal ideal as follows:

Applying the evaluation map twice yields a map $M\otimes_A \Hom(M,
N) \otimes_A N^\vee \to A$, which induces
$$
\epsilon:M\otimes_AN^\vee \longrightarrow E=\Hom(M,N)^\vee
\subseteq S=\S(E).
$$
We can now define the ``generic" determinatal ideal $I_{t+1}
\subseteq S$ to be the ideal generated by the image of the map
$$
\ba{c}
\Lambda^{t+1} M\otimes_A \Lambda^{t+1} (N^\vee) \longrightarrow S\\
m_{0}\wedge\ldots\wedge m_{t}\otimes \varphi_{0}\wedge\ldots\wedge
\varphi_{t}\longmapsto \det(\epsilon(m_j \otimes \varphi_{i}))\in
S.\ea
$$
It is clear that \bnum[(A)] \item the construction yields the
previously defined determinantal ideal (see 4.1) when $N$ is free;
in particular $I_{t+1}$ is the ideal of $(t+1)$-minors of an $n
\times m$ matrix of variables if both $M$ and $N$ are free of
ranks $m$ and $n$, respectively; \item applying the map
$\varphi_f: S \to A$ induced by $f$ (see (3.1), one obtains
$\varphi_f(I_{t+1})=I_{t+1}(f)$, the determinantal ideal of $f$.
\enum

Also notice that $I_{t+1}\subseteq S$ is a homogeneous ideal
generated by forms of positive degree. From item (A) and Lemma
\ref{detideals} below we see that $\dim S_\fq/I_\fq = r-\tau$ with
$r:=nm=\rank E_\fq$. As $I_{t+1}(f)=\varphi_f(I_{t+1})$ by item
(B), the present theorem now follows from Theorem \ref{1.2} in
case (a) and Theorem \ref{1.2n} in case (b). \eproof

\blem\label{detideals} Let $S=A[T_{ij}]$ be the polynomial ring
over a noetherian ring in the variables $T_{ij}$, $1\le i\le n$,
$1\le j\le m$, and let $I_{t+1}$ be the ideal of $(t+1)$-minors of
the generic matrix $T=(T_{ij})$, where $t\le\min\{n,m\}$. Then
$S/I_{t+1}$ is flat over $A$ with geometrically normal fibers of
dimension $nm-\tau $, where $\tau:=(n-t)(m-t)$.

Likewise this holds for the generic symmetric $m\times m$-matrix,
in which case the fiber dimension is
$$
\binom{m+1}{2}-\tau \mbox{ \, with }\tau:=\binom{m-t+1}{2},
$$
and for the generic alternating $m\times m$-matrix, for which the
fiber dimension is
$$
\binom{m}{2}-\tau \mbox{ \, with }t:=2s \mbox{ and }
\tau:=\binom{m-2s}{2}.
$$
\elem

For the proofs we refer the reader to \cite{EN} in the generic
case, \cite[2.1 and 2.4]{Jo} in the symmetric case and \cite[2.1
and 2.5]{JP} in the alternating case.

\bthm \label{gencon} Let $(A,\fm)$ be a noetherian local ring and
let $f:M\to N$ be a linear map of finite $A$-modules with an
isolated singularity of ranks $m$ and $n$, respectively. Assume
that one of the conditions $(a)$, $(b)$ in Theorem
$\ref{gencodim}$ is satisfied. If $\hat A$ is connected in
dimension $d$ and $t\le\min\{n,m\}$, then with $\tau:=(n-t)(m-t)$
the degeneracy locus $D_t(f)$ is connected in dimension $d-\tau$.
\ethm

\bproof  As in the proof of Theorem \ref{gencodim} we let $E$ be
the dual of $\Hom(M,N)$ and $I=I_{t+1}\subseteq \S(E)$ the
``generic" determinantal ideal so that $\varphi_f(I)=I_{t+1}(f)$
in the notation of Theorems \ref{1.6a} and \ref{1.6n}. The map
$V(I)\to \Spec A$ is flat over points $\fp\ne \fm$, and the fibers
are generic determinantal varieties. According to Lemma
\ref{detideals} these fibers are geometrically normal of constant
dimension $nm-\tau$. Now the result follows from Theorems
\ref{1.6a} and \ref{1.6n}. \eproof

Order ideals are a special case of determinantal ideals, namely
the case where $M=A$ and $t=0$. Applying Theorems \ref{gencodim}
and \ref{gencon} to this case we obtain for instance the following
result:

\bcor \label{order} Let $(A,\fm)$ be a noetherian local ring and
let $N$ be a finite $A$-module with an isolated singularity of
rank $n$. Assume that $f \in N$ satisfies one the the following
two conditions $:$ \bnum[$(a)$] \item $f\in \fm N$. \item $A$ has
prime characteristic $p>0$ and the pair $(Af,N)$ is $p$-ample.
\enum Then the following hold $:$ \bnum[$(1)$] \item $\dim
(A/N^\vee(f))\ge \dim A-n$. \item If $\hat A$ is connected in
dimension $d$ then $A/N^\vee(f)$ is connected in dimension $d-n$.
\enum \ecor

\brem (1) One should expect that Theorems \ref{gencodim} and
\ref{gencon} remain valid in arbitrary characteristic if the pair
$(Af, \Hom(M,N))$ is assumed to be ample. A positive answer to
this problem would yield the corresponding projective results (see
Remark \ref{last}(1)) in full generality.

(2) We do not know whether in Theorems \ref{gencodim} and
\ref{gencon} one can also obtain the refined bounds of Bruns
\cite{Br} for determinantal ideals $I_{t+1}(f)$ when
$I_{t+2}(f)=0$. We note that the analogous question for
determinantal ideal of symmetric maps is not true in positive
characteristic as is seen by the simple example of the symmetric
matrix
$$
F=\left(
\ba{cccc}
0&X_1&X_2\\
X_1& 0 & X_3\\
X_2 & X_3 & 0
\ea
\right)
$$
in characteristic 2, see \cite{EHU4}. Therefore one cannot expect
that the solution to this problem follows from the results of
Section 3. \erem

\subsection{Symmetric maps} Let $M$ and $L$ be finite
modules over the noetherian local ring $(A,\fm)$. In this
subsection we will study degeneracy loci and determinantal ideals
of symmetric homomorphisms $f:M\to \Hom_A(M,L)$. By this we mean
that the corresponding map in $\Hom_A(M\otimes_AM,L)$ factors
through the symmetric power $\S^2(M)$.
By $\Hom^s(M,\Hom(M,L))$ we denote the $A$-module of all such
symmetric maps so that 
$$
\Hom^s(M,\Hom(M,L))\cong \Hom(\S^2 M,L)
$$
With this notation we have the following
result:

\bthm \label{symcodim} Assume there exists a prime ideal
$\fq\subseteq A$ with $\dim A=\dim A/\fq$ such that $M_\fq$ and
$L_\fq$ are free of ranks $m$ and $1$, respectively. If  $t\le m$
then with $\tau:=\binom{m-t+1}{2}$ we have $\dim D_t(f)\ge \dim
A-\tau$ in each of the following cases\, $:$ \bnum[$(a)$] \item
$f\in \fm \Hom^s(M,\Hom(M,L))$. \item $A$ has prime characteristic
$p$ and the pair $(Af,\Hom^s(M,\Hom(M,L)))$ is $p$-ample. \enum
\ethm

\bproof With $E:=\Hom^s(M,\Hom(M,L))^\vee$ the morphism $f$ can be
considered as an element in $E^\vee$. As in the proof of Theorem
\ref{gencodim} the idea is to apply Theorems \ref{1.2} and
\ref{1.2n} to the ``generic" symmetric determinantal ideal
$I_{t+1}\subseteq S:=\S(E)$, which can be defined as follows:

Applying the evaluation map twice yields a map $$M\otimes_A
\Hom^s(M, \Hom(M,L)) \otimes_A \Hom(M,L)^\vee \longrightarrow A,$$
which in turn induces
$$
\epsilon: M\otimes_A \Hom(M,L)^\vee \longrightarrow E=\Hom^s(M,
\Hom(M,L))^\vee \subseteq S=\S(E).
$$
We then define the ``generic" symmetric determinatal ideal
$I_{t+1} \subseteq S$ to be the ideal generated by the image of
the map
$$
\ba{c}
\Lambda^{t+1} M\otimes_A \Lambda^{t+1} (\Hom(M,L)^\vee) \longrightarrow S\\
m_{0} \wedge\ldots\wedge m_{t} \otimes
\varphi_{0}\wedge\ldots\wedge \varphi_{t}\longmapsto
\det(\epsilon(m_j \otimes \varphi_{i}))\in S \,.\ea
$$
One easily sees that \bnum[(A)] \item $I_{t+1}$ is the ideal of
$(t+1)$-minors of a symmetric $m \times m$ matrix of variables if
both $M$ and $L$ are free of ranks $m$ and $1$, respectively;
\item applying the map $\varphi_f: S \to A$ induced by $f$ (see
3.1), one obtains the ideal $\varphi_f(I_{t+1})=I_{t+1}(f)$
defining the degeneracy locus $D_t(f)$. \enum

Furthermore $I_{t+1}\subseteq S$ is a homogeneous ideal generated
by forms of positive degree. From property (A) and Lemma
\ref{detideals} one sees that $\dim S_\fq/I_\fq = r-\tau$ with
$r:=\binom{m+1}{2}=\rank E_\fq$. Since
$I_{t+1}(f)=\varphi_f(I_{t+1})$ by property (B), the theorem
follows once again from Theorem \ref{1.2} in case (a) and Theorem
\ref{1.2n} in case (b). \eproof

\bthm \label{symcon} Assume in Theorem $\ref{symcodim}$ that $M$
and $L$ are modules with an isolated singularity of ranks $m$ and
$1$, respectively, and that one of the conditions $(a)$ or $(b)$
is satisfied. If $\hat A$ is connected in dimension $d$ and $t\leq
m$, then with $\tau:=\binom{m-t+1}{2}$ the degeneracy locus
$D_t(f)$ is connected in dimension $d-\tau$. \ethm

\bproof The proof is analogous to that of Theorem \ref{gencon}.
\eproof

\subsection{Alternating maps} Let $M$ and $L$ be finite modules
over the noetherian local ring $(A,\fm)$. The results of the
preceding section hold likewise for alternating maps $f:M\to
\Hom_A(M,L)$, i.e.\ maps corresponding to elements in
$\Hom_A(\Lambda^2M,L)$.
By $\Hom^a(M,\Hom(M,L))$ we denote the $A$-module of all such
alternating maps. The proofs of the next two results are analogous
to those in the symmetric case and are left to the reader.

\bthm \label{altcodim} Assume there exists a prime ideal
$\fq\subseteq A$ with $\dim A=\dim A/\fq$ such that $M_\fq$ and
$L_\fq$ are free $A_{\fq}$-modules of ranks $m$ and $1$
respectively. With $t:=2s\le m$ and $\tau:=\binom{m-2s}{2}$ we
have $\dim D_t(f)\ge\dim A-\tau$ in each of the following cases\,
$:$. \bnum[$(a)$] \item  $f\in \fm \,\Hom^a(M, \Hom(M,L))$. \item
$A$ has prime characteristic $p$ and the pair
$(Af,\Hom^a(M,\Hom(M,L)))$ is $p$-ample. \enum \ethm

\bthm \label{altcon} Assume in Theorem $\ref{altcodim}$ that $M$
and $L$ are modules with an isolated singularity of ranks $m$ and
$1$, respectively, and that one of the conditions $(a)$ or $(b)$
is satisfied. If $\hat A$ is connected in dimension $d$ and $t:=
2s\leq m$, then with $\tau:=\binom{m-2s}{2}$ the degeneracy locus
$D_t(f)$ is connected in dimension $d-\tau$. \ethm

\subsection{Further results}
One can generalize Theorems \ref{1.6}, \ref{1.6a}, \ref{1.6n},
\ref{gencon}, \ref{symcon} and \ref{altcon} to the case where the
modules do not have an isolated singularity. As a sample we state
the following result:

\bcor Let $(A,\fm)$ be a noetherian local ring, $\fa \subseteq A$
an ideal, and $f:M\to N$ be a map of finite $A$-modules that are
locally free on $\Spec A\backslash V(\fa)$ of constant ranks $m$
and $n$, respectively. Assume that with $t\le\min\{m,n\}$ and
$\tau:=(m-t)(n-t)$ we have $\dim A/\fa<\dim A-\tau$. If one of the
conditions $(a)$ or $(b)$ in Theorem $\ref{gencodim}$ is satisfied
and $\hat A$ is connected in dimension $d$, then $D_t(f)\backslash
V(\fa)$ is connected in dimension $d-1-\tau$. \ecor

\bproof We may assume that $A$ is complete, and proceed by
induction on $\dim A/\fa$. In case $\dim A/\fa=0$ this is just
Theorem \ref{gencon}. Now let $\dim A/\fa>0$. By Theorem
\ref{gencodim}  $D_t(f)\backslash V(\fa)$ is not empty and hence
has dimension $>d-1-\tau$. Thus if $D_t(f)\backslash V(\fa)$ is
not connected in dimension $d-1-\tau$, then there exist proper
closed subsets $V_1$, $V_2$ of $D_t(f)\backslash V(\fa)$ with
$D_t(f)\backslash V(\fa)=V_1\cup V_2$ and $\dim V_1\cap V_2 <
d-1-\tau$. For a sufficiently general element $x\in \fm$ we have
$$
\dim A/(xA+\fa)= \dim A/\fa-1<\dim A/xA-\tau
$$
and $\dim V_1\cap V_2\cap V(x) <d-2-\tau$. Moreover, since every
irreducible component of $D_t(f)\backslash V(\fa)$ has dimension
$\geq d -\tau$, the subsets $V_i\cap V(x)$ of $D_t(f)\cap
V(x)\backslash V(\fa)$ are proper. Hence $D_t(f)\cap
V(x)\backslash V(\fa)$ is not connected in dimension $d-2-\tau$.
On the other hand, $A/xA$ is connected in dimension $d-1$ by
Grothendieck's connectedness theorem, and conditions (a) and (b)
pass from $A$ to $A/xA$ according to Remark \ref{rempample}(4).
Therefore the induction hypothesis shows that $D_t(f)\cap
V(x)\backslash V(\fa)$ is connected in dimension $d-2-\tau$, which
is a contradiction. \eproof

\brem Similar results hold for symmetric and alternating maps. We
leave the straightforward formulations and proofs to the reader.
\erem

One can also generalize the preceding results in the spirit of
Theorems \ref{1.2}, \ref {1.6} and Corollary \ref{1.8}. Again, we
only give a sample of such a generalization. The special case of
order ideals had been treated in \cite{EHU3}.

\bcor\label{generalization}  Theorem $\ref{gencodim}$ remains
valid if we replace condition $(a)$ by the weaker condition that
there exists a local map $h:A\to A'$ to a noetherian local ring
$(A',\fm')$ so that $\fq=h^{-1}(\fq')$ for some prime $\fq'$ of
$A'$ and $$f\otimes 1\in \fm'(\Hom_A(M,N)\otimes_AA'/{\rm
torsion}).$$\ecor

This is easily derived as in the proof of Theorem \ref{gencodim}
using Corollary \ref{1.8}. Again we only mention the fact that
similar generalizations hold for symmetric and alternating maps;
the straightforward formulations and proofs are left to the
reader.

Finally we apply these results to maps of vector bundles on
projective schemes.

\bcor \label{genproj} Let $X$ be a projective $K$-scheme over
field $K$ of characteristic $p>0$ and let $f:\cM\to \cN$ be an
$\cO_X$-linear map of locally free sheaves of ranks $m$ and $n$,
respectively. If \, $\cHom(\cM,\cN)$ is cohomologically $p$-ample
and $t\le\min\{m,n\}$, then with $\tau:=(m-t)(n-t)$ the following
hold $:$ 
\bnum[$(1)$] 
\item If $\dim X\ge\tau$ then $D_t(f)\ne
\emptyset$. 
\item If $X$ is connected in dimension $d$ then
$D_t(f)$ is connected in dimension $d-\tau $. 
\enum 
\ecor

\bproof 
Write $X=\Proj R$, where $R=K[R_1]$ is a finitely
generated graded algebra over $K=R_0$. The $R$-modules
$$
M=\bigoplus_{n\ge 0} H^0(X, \cM(n))\quad \mbox{and}\quad
N=\bigoplus_{n\ge 0} H^0(X,\cN(n))
$$
are then finitely generated
over $R$ and $\cM=\tilde M$,
$\cN=\tilde N$ are the sheaves associated to $M$, $N$,
respectively. As $\cM$, $\cN$ are vector bundles on $X$, the
modules $M$, $N$ have an isolated singularity. Let
$I_{t+1}(f_0)\subseteq R$ be the $(t+1)$st determinantal ideal of
the map $f_{0}:M\to N$ induced by $f$. Clearly $f_0$ is an element
of degree 0 in $\Hom_R(M,N)$ and so by Proposition
\ref{pample}(1) and Remark \ref{rempample}(2) the pair $(R f_0,
\Hom_{R}(M,N))$ is
$p$-ample. Let $\fm$ be the homogeneous maximal ideal of $R$.
Applying Theorem \ref{gencodim} to the local ring $A=R_{\fm}$ and
the localized map $f_0:M_\fm\to N_\fm$, we obtain that $\dim
A/I_{t+1}(f_{0})A\ge \dim A-\tau$ and so $\dim R/I_{t+1}(f_0)R\ge
\dim R-\tau$. As $\cI_{t+1}(f)\subseteq \cO_{X}$ is the ideal
sheaf induced by $I_{t+1}(f_{0})$, the first part follows.
Likewise part (2) is an easy consequence of Theorem \ref{gencon}.
Here one also uses the fact that a projective scheme is connected
in dimension $d$ if and only if the local ring at the vertex of
the affine cone is connected in dimension $d+1$ (see \cite[Lemma
3.1.9]{FOV}). \eproof

In the case of symmetric or alternating maps one can apply
Theorems \ref{symcodim}, \ref{symcon}, \ref{altcodim} and
\ref{altcon} instead. This leads to the next result. Here ${\bf
D}^n$ denotes the $n$th divided power functor.

\bcor 
\label{alt-symproj} 
Let $X$ is a projective $K$-scheme over
a field $K$ of characteristic $p>0$. Assume that $\cM$ is a
locally free $\cO_X$-module of rank $m$ and $\cL$ is a line bundle
on $X$. Then for $t\le m$ the following hold $:$ \bnum[$(1)$]
\item If $f:\cM^\vee\to  \cM\otimes_{\cO_X}\cL$ is a symmetric map
and \, ${\bf D}^2\cM\otimes_{\cO_X}\cL$ is cohomologically
$p$-ample, then statements $(1)$ and $(2)$ of Corollary
$\ref{genproj}$ hold with $\tau:=\binom{m-t+1}{2}$. 
\item If $f
:\cM^\vee\to \cM\otimes_{\cO_X}\cL$ is an alternating map and
$\Lambda^2\cM\otimes_{\cO_X}\cL$ is cohomologically $p$-ample,
then statements $(1)$ and $(2)$ of Corollary $\ref{genproj}$ hold
with $t:=2s$ and $\tau:=\binom{m-2s}{2}$.
\enum 
\ecor

\brem\label{last} (1) Suppose that $X$ is a projective variety
over the field $K$. If $K$ has characteristic $0$ or $X$ is
smooth, the conclusion of Corollary \ref{genproj} has been shown
by Fulton and Lazarsfeld \cite{FL} under the assumption that
$\cHom(\cM,\cN)$ is ample. Similarly, Corollary \ref{alt-symproj}
with the assumption of $p$-ampleness replaced by ampleness is due
to Tu \cite{Tu} provided that $K$ has characteristic $0$ and the
maps have even ranks.

(2) Using Theorems \ref{gencodim} and \ref{gencon} we also obtain
an analogue of Corollary \ref{genproj} in arbitrary characteristic
if we replace the condition on the cohomological $p$-ampleness of
the vector bundle $\cH:=\cHom (\cM,\cN)$  by the assumption that
the element $f\in H^0(X, \cH)$ can be written as a linear
combination of sections in $H^0(X, \cH(-1))$ with coefficients in
$H^0(X, \cO_X(1))$. A similar remark applies to Corollary
\ref{alt-symproj}. 
\erem

\end{document}